\documentclass[12pt, a4paper]{article}
\usepackage{amsmath}
\usepackage{amsthm}
\usepackage{amsfonts}
\usepackage{todonotes}
\usepackage{nicefrac}
\usepackage{dsfont}
\usepackage{amssymb}
\usepackage{latexsym}
\usepackage{mathtools}
\usepackage{verbatim}
\usepackage{enumitem}
\usepackage{rotating}
\usepackage{subcaption}
\usepackage{float}
\usepackage{multirow}
\usepackage{stmaryrd}
\usepackage{tikz}
\usepackage{upgreek}

\usepackage[backend=bibtex,url=false,isbn=false,style=alphabetic,maxnames=5,firstinits=true]{biblatex}
\renewbibmacro{in:}{}
\AtEveryBibitem{%
  \clearlist{language}%
  \clearfield{issn}%
  \clearfield{month}%
}
\DeclareNameAlias{default}{first-last}

\let\oldcite\cite
\renewcommand{\cite}[2][\empty]{\oldcite{#2}}

\DefineBibliographyStrings{english}{%
  bibliography = {References},
}

\usepackage[
    pdftex,
		bookmarks,
		bookmarksopen,
		bookmarksopenlevel=0,
		bookmarksnumbered,
		breaklinks,
		colorlinks=true,
		linkcolor=linkcolor,
		citecolor=citecolor,
		urlcolor=urlcolor,
	]{hyperref}

\hypersetup{
	pdftitle    = {Dissertation},
	pdfauthor   = {Andreas Vollmer},
	pdfsubject  = {dissertation},
	pdfcreator  = {\LaTeX{} with `book' class},
	pdfproducer = {pdf\LaTeX},
	pdfkeywords = {first integrals, Killing tensors, stationary axisymmetric vacuum metrics, sub-Riemannian structure}
}

\definecolor{linkcolor}{rgb}{0.5,0.0,0.0}
\definecolor{citecolor}{rgb}{0.0,0.5,0.0}
\definecolor{urlcolor} {rgb}{0.0,0.0,0.5}

\theoremstyle{plain} 
\newtheorem{theorem}{Theorem}

\newtheorem*{result*}{Result}

\newtheorem{algorithm}{Algorithm}

\newcommand{\circled}[1]{\ensuremath{\text{($#1$)}}}
\newcommand{\rg}{\mathrm{rk}}

\newcommand{\Ric}{\ensuremath{\mathrm{Ric}}}

\newcommand{\TaM}{\ensuremath{T^\ast\!M}}

\newcommand{\E}{\ensuremath{\mathcal{E}}}  
\newcommand{\Seqn}{\ensuremath{\mathcal{S}}} 
\newcommand{\Sd}{\ensuremath{\mathcal{S}_d}} 
\newcommand{\meqns}[2]{\ensuremath{m_{#1,#2}}}  
\newcommand{\nvars}[2]{\ensuremath{n_{#1,#2}}}  

\newcommand{\parity}{\ensuremath{\wp}}

\newcommand{\trivials}[1]{\ensuremath{\Lambda_{#1}^0}}
\newcommand{\upperbound}[2][\empty]{\ifthenelse{\equal{#1}{\empty}}{\ensuremath{\overline{\Lambda}_{#2}}}{\ensuremath{\overline{\Lambda}_{#2}^{(#1)}}}}

\newcommand{\intrange}[1]{\llbracket{#1}\rrbracket}

\newcommand{\wrt}{w.r.t.}

\newcommand{\ie}{i.e.}
\newcommand{\eg}{e.g.}
\newcommand{\cf}{cf.}

\newcommand{\RR}{\ensuremath{\mathds{R}}}

\begin{document}

\title{Killing tensors in stationary and axially symmetric space-times}
\author{Andreas Vollmer (andreas.vollmer@uni-jena.de)}
\maketitle

\section*{Abstract}
We discuss the existence of Killing tensors for certain (physically motivated) stationary and axially symmetric vacuum space-times.
We show nonexistence of a nontrivial Killing tensor for a Tomimatsu-Sato metric (up to valence~7), for a C-metric (up to valence~9) and for a Zipoy-Voorhees metric (up to valence~11).

The results are obtained by mathematically completely rigorous, nontrivial computer algebra computations with a huge number of equations involved in the problem.

\section{Introduction}\label{sec:introduction}

Let $(M,g)$ be a 4-dimensional manifold with Lorentzian metric~$g$ of signature~(+,+,+,--).
A \emph{Killing tensor} of valence~$d$ on~$M$ is a symmetric tensor field $K$ whose symmetrized covariant derivative vanishes,
\begin{equation}\label{eqn:killing}
 \nabla_{(j}K_{i_1\dots i_d)}=0.
\end{equation}
Here, $\nabla$ denotes the Levi-Civita connection of~$g$ and the components $K_{i_1\dots i_d}$ smoothly depend on the position coordinates.

The metric~$g$ itself trivially is a Killing tensor of valence~2.
Killing vectors with lowered indices are valence-1 Killing tensors.
Given two Killing tensors, we can construct a new Killing tensor by forming their symmetrized product.
Since Equation~\eqref{eqn:killing} is linear in the components of $K$, the linear combination of Killing tensors of the same valence is also a Killing tensor.

Killing tensors correspond to first integrals of the geodesic flow that are homogeneous polynomials in the momenta:
An integral of the geodesic flow is a function $I:\TaM\to\RR$ such that its Poisson bracket with the Hamiltonian $H:\TaM\to\RR$, $H=g^{ij}p_ip_j$, vanishes, \ie\
\begin{equation}\label{eqn:poisson}
 \{H,I\}=\sum_i \frac{\partial H}{\partial x^i}\frac{\partial I}{\partial p_i}-\frac{\partial H}{\partial p_i}\frac{\partial I}{\partial x^i} = 0.
\end{equation}
It is well known that for a Killing tensor $K$ the function,
\begin{equation}\label{eqn:integral}
 I_K(x,p) = \sum K^{i_1\dots i_d}(x)\,p_{i_1}\cdots p_{i_d},
\end{equation}
is an integral (here the $K^{i_1\dots i_d}$ denote components of a Killing tensor with raised indices).

Two Killing tensors are in involution if their corresponding integrals commute \wrt\ the standard Poisson bracket on \TaM.

The requirement \eqref{eqn:integral} is equivalent to a system of partial differential equations (PDEs) on the coefficients of~$I$, \ie\ on the $K^{i_1\dots i_d}$.
These equations are coefficients of~$\{I,H\}=0$ and therefore obviously polynomial in the momenta.

Killing tensors appear in many contexts, \eg\ in mechanics, mathematical relativity, integrability or differential geometry.
In general relativity, the geodesic flow determines trajectories of free-falling particles. Integrals of the geodesic flow provide \emph{constants of the motion} of such particles.
Hamiltonian systems on 4-dimensional manifolds with four functionally independent integrals in involution are Liouville integrable. For such systems the orbits are restricted to tori (compact case) or cylinders, and the equations of motion can be solved by quadrature.

It has recently been shown by Kruglikov \& Matveev that a generic metric does not admit non-trivial Killing tensors~\cite{kruglikov_generic_2015}.
However, many examples in physics do admit nontrivial Killing tensors. Such examples are in some sense the most important ones, \eg\ the classical Kepler problem.
In the context of stationary and axially symmetric vacuum metrics, the most prominent example is the family of Kerr metrics. These metrics are used as a model for the space-time around rotating neutron stars and black holes.
In particular, the Schwarzschild metric is the static limit of the Kerr family.

A metric satisfies the \emph{vacuum condition} if it is Ricci-flat, \ie\ if its Ricci tensor vanishes.
This requirement is a system of partial differential equations on the components of the metric.
For our examples, this requirement is automatically satisfied. See \cite{vollmer_reducibility_2015} for a result on Killing tensors of valence~3 on arbitrary static and axially symmetric vacuum space-times. This reference makes explicit use of the vacuum condition.
Physically, the vacuum condition $\Ric(g)=0$ is a fair assumption in the exterior region of stationary and axially symmetric astrophysical objects when we ignore electromagnetic fields.

The present paper generally follows the method used by Kruglikov \& Matveev in \cite{kruglikov_nonexistence_2011}, which is based on an algorithm. We improve this algorithm and achieve considerably higher computational efficiency.
The method can prove non-existence of nontrivial Killing tensors. In case there are additional Killing tensors, the method finds them.

In our computations, we have to deal with up to more than 10,000 equations and unknowns, which necessitates the use of computer algebra.
This has also been done in \cite{kruglikov_nonexistence_2011} where nonexistence of a nontrivial Killing tensor of valence up to 6 is proven for the Darmois metric (in \cite{kruglikov_nonexistence_2011} it is called the Zipoy-Voorhees metric with $\delta=2$).
The reachable valence of the Killing tensors is, however, only restricted by computer strength.
Our code is based on \cite{kruglikov_nonexistence_2011}, but we employ additional tricks and achieve for the same metric a nonexistence proof up to valence 11. Note that our result on the Darmois metric does not follow from \cite{maciejewski_nonexistence_2013} where only analytic integrals are considered. See also Section~\ref{sec:darmois}.

Moreover, we implement the method for the first time for a non-static metric. Specifically, we prove for a certain Tomimatsu-Sato metric that there are no additional Killing tensors up to valence~7.

The paper is organized as follows.
In Section~\ref{sec:methods}, we give a description of the method we use. The algorithm is summarized on page~\pageref{alg:sav}.
%
Then, in Section~\ref{sec:examples}, examples are given to exemplify application of the method.
We investigate a (non-static) Tomimatsu-Sato metric in Theorem~\ref{thm:TS-7}, a Zipoy-Voorhees metric (the Darmois solution) in Theorem~\ref{thm:ZV-11} and a C-metric in Theorem~\ref{thm:C-metric}. As an example for a case of existence of a nontrivial Killing tensor, we discuss the Kerr metric in Section~\ref{sec:kerr}.

\section{Method}\label{sec:methods}
The general idea is classical and sometimes called Cartan-K\"ahler method or prolongation-projection method.
It allows us to deal with a nontrivial overdetermined system of PDEs, say $\Seqn$, by straightforward, but computationally challenging, algebraic calculations.
Actually, the system $\Seqn$ is linear in our cases, and therefore our computations are also linear-algebraic.

The basic procedure is as follows. Consider the differential consequences of~$\Seqn$, \ie\ differentiate the equations \wrt\ the independent coordinate variables. The system of equations resulting after $k$ differentiations is called the $k$-th \textit{prolongation}. An overdetermined system $\Seqn$ of PDEs is called \emph{finite} or \emph{of finite type} if highest derivatives of the unknowns can be expressed through lower derivatives after a finite number of differentiations.
The existence of Killing tensors is equivalent to the existence of solutions for an overdetermined system of PDEs of finite type.
Consider derivatives of the unknown functions (in our case the coefficients in~$I$ and their derivatives) as new, independent \textit{unknowns}. The equations are then algebraic equations on the unknowns. This algebraic system admits at least the solutions corresponding to solutions of the initial differential problem. Solving the algebraic problem therefore leads to an upper bound to the number of Killing tensors. Nonexistence of nontrivial Killing tensors is proven if this upper bound coincides with the number of trivial Killing tensors.


In general, stationary and axially symmetric vacuum metrics can be written in so-called Lewis-Papapetrou coordinates $(x,y,\phi,t)$,
\begin{equation}\label{eqn:metric-sav}
 g_\text{SAV} = e^{2U}\,\left( e^{-2\gamma}\,\left(dx^2+dy^2\right) +x^2\,d\phi^2 \right) +e^{-2U}\,(dt+A\,d\phi)^2,
\end{equation}
with three parametrizing functions $U(x,y)$, $\gamma(x,y)$ and $A(x,y)$.
For the metrics in question, these parametrizing functions are fixed.
In our examples, it is actually convenient not to use the coordinates of~\eqref{eqn:metric-sav}. Instead, we use slightly different coordinate choices, see Equations~\eqref{eqn:tomimatsu-sato}, \eqref{eqn:zipoy-voorhees} and \eqref{eqn:c-metric}.

All our coordinate choices have the following property:
The coordinates are adapted to the obvious symmetries of the metric.
We clearly see that the metrics are invariant under rotations $\phi\mapsto\phi+\phi_0$ and time translations $t\mapsto t+t_0$.
These symmetries correspond to the Killing vectors $\partial_t$ and $\partial_\phi$.

Consider level surfaces of the integrals $p_t$ and $p_\phi$ corresponding to the two Killing vectors $\partial_t$ and $\partial_\phi$, respectively. For regular values, these level surfaces are submanifolds.
They are endowed with the coadjoint action of the symmetry group, and the quotient space $M_\text{red}$ is known as the symplectic quotient.
It inherits a natural symplectic form and a Hamiltonian from the original manifold \cite{whittaker_treatise_1904,marsden_lectures_1992}.
The Hamiltonian on this \emph{reduced} space is, however, inhomogeneous, and therefore the initial problem of finding Killing tensors turns into the problem of finding \emph{inhomogeneous integrals} of the geodesic flow for a reduced Hamiltonian of the form $H=K_g+V$. The first term $K_g$ is called the \emph{kinetic term} and corresponds to a Killing tensor on $M_\text{red}$. The function $V:M\to\mathds{R}$ is called the \emph{potential}. The coordinates $x,y$ (and the respective momenta) are called \emph{non-ignorable}, while $\phi,t$ ($p_\phi,p_t$) are called \emph{ignorable} \cite{whittaker_treatise_1904,carter_hamilton-jacobi_1968}.

We are interested in functions $I:\TaM_\text{red}\to\mathds{R}$ such that the Poisson equation~\eqref{eqn:poisson} for the Hamiltonian flow defined by $H$ is satisfied, \ie\ such that $\{H,I\}=0$.
This equation is an inhomogeneous polynomial in the momenta $(p_x,p_y)$.
We introduce the following notion of parity for such polynomials:
We say that the polynomial is of even (odd) parity if all its homogeneous components have even (odd) degree in the momenta $(p_x,p_y)$.

Since~$H$ is even in the momenta $(p_x,p_y)$, we can consider integrals~$I$ of odd and even parity in these momenta separately, \cf\ \cite{kruglikov_nonexistence_2011, hietarinta_direct_1987}.
Considering homogeneous components of~\eqref{eqn:poisson}, one obtains a list of polynomial equations
\begin{subequations}\label{eqns:decomposition-ni-degree}
\begin{align}
 \E_0\,	& & \{T,I^\circled{d}\}	&= 0 \\
 \E_1\,	& & \{T,I^\circled{d-2}\} +\{V,I^\circled{d}\}	&= 0 \\
 \vdots\,\,\,	& & 										&\hspace{-0.05cm}\dots   \nonumber   \\
 \smash{\raisebox{-0.3cm}{$\E_\text{fin}$}}	& & 
 			\smash{\begin{array}{l}
 				\null\\[\jot]
 				\left\{\begin{array}{l}
    			\{T,I^\circled{0}\} +\{V,I^\circled{2}\} \\[\jot]
    			\{V,I^\circled{1}\}
   			\end{array}\right.
   			\end{array} }
   										&= 0 \quad\text{(even parity branch)} \\
    	& &							&= 0 \quad\text{(odd parity branch)}
\end{align}
\end{subequations}
Here, we denote by $I^\circled{r}$ the homogeneous polynomial component of~$I$ of degree~$r$ \wrt\ the momenta~$(p_x,p_y)$.

The polynomials $\E_k$ can be further decomposed \wrt\ the \emph{ignorable} momenta $p_\phi$ and $p_t$. We have two ignorable momenta, and thus we can decompose the $k$-th equation $\E_{k-1}$ into $2k-1$ new polynomial equations.
We can consider every polynomial equation $\E_l$ as corresponding to a collection of equations provided by the coefficients \wrt\ $p_\phi$ and~$p_t$. Denote these equations by $\E_l$ as well.

We follow the Cartan-K\"ahler method as explained above. We differentiate each $\E_l$ \wrt\ $(x,y)$ and consider derivatives of the components $K^{i_1\dots i_d}$ of~$I$ (\ie\ the unknown functions) as new, independent unknowns. Since the metric is given explicitly in our examples, we obtain a linear-algebraic problem. The number of solutions of the linear problem can be determined by computing the rank of a (huge) matrix.
Solutions of system of PDEs corresponding to~\eqref{eqn:integral} are equivalent to solutions of the linear problem. Therefore the matrix rank is an upper bound to the number of Killing tensors that the metric admits.
Since the matrix dimensions are huge in our cases, we need to find tricks that can speed up the computations, particularly the required rank computation.

For the metrics in question, highest derivatives of the unknowns can be expressed through lower derivatives after differentiating the equations~$\E_l$ for $d+1$ times~\cite{wolf_structural_1997}. We need $d$ differentiations to find the lowest possible upper bound~\cite{kruglikov_nonexistence_2011, vollmer_thesis_2016}.
Specifically, write the Poisson bracket as the polynomial
\begin{equation}\label{eqn:poisson-for-sav}
  \{H,I\}=\sum_{i=0}^{d+1} \sum_{j=0}^i \sum_{k=0}^{d+1-i}
              P^\circled{i,j}_k\ p_1^{i-j}\,p_2^j\,p_3^k\,p_4^{d+1-i-k},
\end{equation}
where each \smash{$P^\circled{i,j}_k=0$} represents an equation in the system of PDEs.
The prolongated system is obtained by taking derivatives of the \smash{$P^\circled{i,j}_k$}.
For given $i$, $j$ and $k$, denote the resulting equation obtained after~$m$ differentiations (with $\mu$ derivatives \wrt\ $x$ and $m-\mu$ \wrt\ $y$) by $P^\circled{i,j,m,\mu}_k$.
We denote by $\intrange{a,b}$ the integers between (and including)~$a$ and~$b$.
The indices run over the following values:
$i\in\intrange{0,d+1},\ j\in\intrange{0,i},\ k\in\intrange{0,d+1-i}$, and $\mu\in\intrange{0,m}$.

For integrals of pure parity \wrt\ $(p_x,p_y)$, many $P^\circled{i,j}_k$ are zero.
Particularly, if we consider integrals of odd (even) $(p_1,p_2)$-parity, then only \smash{$P^\circled{i,j}_k$} with even (odd) value of $i$ can be non-zero.
Now, the unknown functions are the coefficients in the polynomial that represents~$I$,
\begin{equation}
  I = \sum_{\mathclap{\substack{i=0\\ \parity(i)=e}}}^d\ \sum_{j=0}^i\ \sum_{k=0}^{d-i}\ I^\circled{i,j}_k\ p_1^{i-j}\,p_2^j\,p_3^k\,p_4^{d-i-k}
  \qquad\text{with}\ e=0\text{ or }e=1
\end{equation}
and where $\parity(i)$ denotes the parity of the integer $i$.
For the derivatives of the unknown functions, we use a notation analogous to that for the \smash{$P^\circled{i,j}_k$}, namely $I^\circled{i,j,m,\mu}_k$, with $i\in\intrange{0,d},\ j\in\intrange{0,i},\ k\in\intrange{0,d-i},\ \mu\in\intrange{0,m}$.
Here, $m$ denotes the order of partial differentiation, and $\mu$ is the order of differentiation \wrt\ the coordinate $x$.

\paragraph{Structuring the equations and unknowns.}
With the considerations just made, we can organize the equations and unknowns into a tabular structure.
For the equations, consider all sets~$\E_l$ and denote them in one column with~$l$ indexing the rows. Then, put their differential consequences in columns to the right, \ie\ the first derivatives of~$\E_0$ \wrt\ $(x,y)$ are in the first row of the second column and so forth (\cf\ Figure~\ref{fig:ordering-eqns}).

\begin{figure}
 \resizebox{\textwidth}{!}{  \begin{tikzpicture}[x=1.5cm,y=-0.8cm]
    \tikzstyle{every node}=[font=\small]
    \draw[thin,black] (0,0) rectangle (10,7);
    \draw[thin,black] (0,0)--(1,1);
    \node at (0.7,0.3) {\small $m$};
    \node at (0.3,0.7) {\small $l$};
    \draw[thin,black] (0,1)--(10,1);
    \draw[thin,black] (0,2)--(10,2);
    \draw[thin,black] (0,3)--(10,3);
    \node at (0.5,4) {\vdots};
    \node[above] at (5.5,4) {\vdots};
    \node at (5.5,4.1) {\vdots};
    \node[left] at (5.5,4) {\dots};
    \node[right] at (5.5,4) {\dots};
    \node at (9.5,4) {\vdots};
    \draw[thin,black] (0,5)--(10,5);
    \draw[thin,black] (0,6)--(10,6);
    \draw[thin,black] (1,0)--(1,7);
    \draw[thin,black] (2,0)--(2,7);
    \draw[thin,black] (3,0)--(3,7);
    \draw[thin,black] (4,0)--(4,7);
    \node at (5.5,0.5) {\dots};
    \node at (5.5,6.5) {\dots};
    \draw[thin,black] (7,0)--(7,7);
    \draw[thin,black] (8,0)--(8,7);
    \draw[thin,black] (9,0)--(9,7);
    \node at (1.5,0.5) {$0$};
    \node at (2.5,0.5) {$1$};
    \node at (7.5,0.5) {$M\!\!-\!\!2$};
	\node at (8.5,0.5) {$M\!\!-\!\!1$};
	\node at (9.5,0.5) {$M$};
	\node at (0.5,1.5) {$0$};
	\node at (0.5,2.5) {$1$};
	\node at (0.5,5.5) {$\hat{L}\!\!-\!\!1$};
	\node at (0.5,6.5) {$\hat{L}$};
	\node at (1.5,1.5) {$\E_0$};
	\node at (1.5,2.5) {$\E_1$};
	\node at (1.5,5.5) {$\E_{\hat{L}\!-\!1}$};
	\node at (1.5,6.5) {$\E_{\hat{L}}$};
	\node at (2.5,1.5) {$\frac{\partial\E_0}{\partial x^J}$};
	\node at (7.5,1.5) {$\frac{\partial^{M\!-\!2}\E_0}{\partial x^J}$};
	\node at (8.5,1.5) {$\frac{\partial^{M\!-\!1}\E_0}{\partial x^J}$};
	\node at (9.5,1.5) {$\frac{\partial^M\E_0}{\partial x^J}$};
  \end{tikzpicture}}
 \caption{Ordering scheme for the equations, with $(x^1,x^2)=(x,y)$.}
 \label{fig:ordering-eqns}
\end{figure}

The unknowns are the coefficient functions of $I$ \wrt\ momenta. They can be organized in a way similar to the equations.
Let $2l=i-d+\tilde{e}$, where $\tilde{e}\in\{0,1\}$ is the parity of $d+e$.
Then, we first arrange the \smash{$I^\circled{i,j,m,\mu}_k$} according to the value of $l$.
For equal values of $l$, we then arrange the unknowns according to the order $m$ of differentiation.
Note that the resulting table for the unknowns has one column more than the table for the equations.
The reason is that the system of PDEs following from~\eqref{eqn:killing} is of first order.

\paragraph{Elimination scheme}
Let us have a closer look at the structure of these tables.
We regard derivatives of the unknown functions as new, independent unknowns.
We observe that in the table of equations the $(l,m)$-cell shares the unknowns $I^\circled{i,j,m+1,\mu}_k$ with the $(l,m+1)$-cell if $i=2l+d-\tilde{e}$.
Together with the structuring of the unknowns obtained above, this pattern suggests to solve the linear system of equations stepwise.
In principle, one can handle one cell of the table of equations at a time, and iteratively replace unknowns.

However, we are not going to follow this prescription entirely.
Instead, only those equations \smash{$P^\circled{i,j,m,\mu}_k$} that are monomial in the respective unknowns \smash{$I^\circled{i,j,m+1,\tilde{\mu}}_{\tilde{k}}$} will be taken into account. Yet, this will be done iteratively, so a maximal number of substitutions can be achieved.
This partial solution of the system reduces the number of equations and unknowns considerably and therefore improves performance in the following steps.

\paragraph{Computing the number of Killing tensors}

We need to identify the number of solutions of the obtained linear-algebraic system. This system is described by a matrix and we have to compute the rank of this matrix.
On a computer, we can do this by choosing a point of reference in which we complete the computations. Actually, a little caution is necessary since we need the number of solutions for the generic matrix system; in non-generic points, the rank of the matrix may drop. We also restrict to rational reference points.
In case the expressions are rational in the coordinates $x$ and $y$, we can rewrite the equations such that the coefficients become integer numbers.

After choosing a reference point, our freedom to add arbitrary multiples of known integrals can be made use of to further eliminate unknowns from the system.
Specifically, we can set all \smash{$I^\circled{i,0,0,0}_k=0$}.

For the remaining matrix problem, we can determine the number of solutions from the dimensionality of the matrix kernel.
We use usual Gau{\ss} elimination for this computation.
Better algorithms might be available for particular situations.

The computation of an upper bound to the number of involutive Killing tensors can now be performed algorithmically:
\begin{algorithm}\label{alg:sav}
\ 
 \begin{enumerate}[label={(\roman*)}, align=right]
  \item\label{step:initial} Consider the two pure-parity integrals \wrt\ $(p_x,p_y)$ separately. Compute the differential consequences of the corresponding differential systems up to $d$-th prolongation. Consider the corresponding algebraic problem.
	\item Choose a generic point $P$ and evaluate the algebraic system at this point. Add multiples of the known integrals to set as many of the unknowns as possible to zero (in $P$).
	\item Perform the elimination scheme as discussed above.
	\item If possible, rewrite the matrix system such that the coefficients are integers. Determine the kernel dimension.
 \end{enumerate}
\end{algorithm}
The algorithm confirms nonexistence of an additional integral if the matrix has full rank.


\section{Examples}\label{sec:examples}

Let us explore some examples with the method we developed in the previous section.
Coordinates are chosen according to the specific problem and are not identical to those in Equation~\eqref{eqn:metric-sav}.
However, they are still adjusted to the symmetries, \ie\ to stationarity and axial symmetry.

The computation times have been achieved on a desktop computer with a 3.4GHz processor and 32GB RAM. The computations were performed using Maple 18.

\subsection{Tomimatsu-Sato metrics}
The Tomimatsu-Sato family generalizes the Kerr metric.
Its static subclass is the Zipoy-Voorhees family, which contains the Schwarzschild metric as its Kerr limit.

\subsubsection{A non-static example}
We begin with a non-static case and consider a Tomimatsu-Sato metric with perturbation parameter $\delta=2$. In the Ernst-Perj\'es representation, it has the general form \cite{ernst_representation_1976,perjes_factor_1989,manko_physical_2012}
\begin{equation}\label{eqn:tomimatsu-sato}
\resizebox{0.9\textwidth}{!}{
  $g_\text{TS} = \kappa^2\,f^{-1}\ \bigg( e^{2\gamma}\,(x^2-y^2)\bigg(\frac{dx^2}{x^2-1}+\frac{dy^2}{1-y^2}\bigg)
  																				+(x^2-1)(1-y^2)\,d\phi^2 \bigg)
  										-f\,(dt-\omega\,d\phi)^2$.
}
\end{equation}
where the functions $f$, $\gamma$ and $\omega$ are defined by
\begin{subequations}
\begin{align}
 f						&= \frac{\mu^2-(x^2-1)(1-y^2)\sigma^2}{\mu^2+\mu\nu-(1-y^2)((x^2-1)\sigma^2-\sigma\tau)}, \\
 e^{2\gamma}	&= \frac{\mu^2-(x^2-1)(1-y^2)\sigma^2}{p^4\,(x^2-y^2)^4}, \\
 \omega				&= -\frac{\kappa\,(1-y^2)((x^2-1)\sigma\nu+\mu\tau)}{\mu^2-(x^2-1)(1-y^2)\sigma^2}, \\
\intertext{and where $\mu$, $\nu$, $\sigma$ and $\tau$ are the polynomials}
 \mu					&= p^2\,(x^2-1)^2+q^2\,(1-y^2)^2, \\
 \nu					&= 4x\,(px^2+2x+p), \\
 \sigma				&= 2pq\,(x^2-y^2), \\
 \tau					&= -4qp^{-1}\,(1-y^2)(px+1).
\end{align}
\end{subequations}
In addition,~$p$ and~$q$ have to obey the restriction, $p^2+q^2=1$.
The other free parameter, $\kappa$, can in principle be removed through redefinition of some quantities, but we keep it in view of \cite{manko_physical_2012}.
We study the particular example with parameter values $\delta=2$, $\kappa=2$, and $p=\nicefrac{3}{5}\ (q=\nicefrac{4}{5})$.
These parameters have also been chosen in \cite{manko_physical_2012}, where some physical properties of the Tomimatsu-Sato metric for $\delta=2$ are discussed.

\begin{theorem}\label{thm:TS-7}
For this Tomimatsu-Sato metric, there is no additional independent Killing tensor of valence~$d\leq7$ that is in involution with the trivial Killing tensors $d\phi$, $dt$, and the metric.
\end{theorem}

Indeed, the following table shows that after $d=7$ steps of prolongation the algorithm yields the upper bound 20 (sum of the \smash{$\upperbound[d]{d}$}, which denote the obtained upper bound for degree~$d$ and after $d$ prolongation steps).
The number $\trivials{d}$ of trivial integrals is given by the formula
\[
 \trivials{d} = \sum_{l=0}^{\lfloor\nicefrac{d}{2}\rfloor} \left(\begin{array}{c} D+d-2l-3 \\ d-2l\end{array}\right)
 \underset{\text{D=4}}{\overset{\text{d=7}}{=}} 20.
\]
Both numbers coincide and this confirms the nonexistence of an additional Killing tensor.

\begin{center}
\textbf{Results Tomimatsu-Sato metric}\\
\textit{$\delta=2$, $\kappa=2$, $p=\nicefrac{3}{5}\ (q=\nicefrac{4}{5})$}
\vskip 0.2cm
\begin{tabular}{|c|c|ccccccc|}
\hline
d                				  & e & \upperbound[d]{d} & $\meqns{d}{d}$ & $\nvars{d}{d}$ & rows of $M$ & columns of $M$ & $\rg(M)$ & time  \\ \hline
\multirow{2}{*}{\large7}	& 0 & 20    & 2880 & 2700 & 556 & 356 & 356 & 21h	   \\
                   				& 1 & 0     & 3060 & 2700 & 776 & 416 & 416 & 24h    \\ \hline
\end{tabular}
\end{center}
Recall that $e$ is the parity of the integral \wrt\ the momenta $(p_x,p_y)$. Thus, for given valence~$d$, we have to compute two separate branches $e=0$ and $e=1$.

By $M$, we denote the matrix obtained after Step~\ref{step:initial} of the algorithm.
The symbols $\meqns{d}{d}$ and $\nvars{d}{d}$ denote, respectively, the number of equations and unknowns of the initial matrix system obtained after Step~\ref{step:initial} of the algorithm, \ie\ for degree $d$ and after $d$ prolongation steps. The point of reference for the computations is $(x,y)=(\nicefrac12,2)$. The last column provides the (approximate) computation times, \cf\ also Section~\ref{par:efficiency}.
\medskip

One might wonder why the number $\nvars{d}{d}$ coincides for both branches. The reason is that for the number of unknowns the following formula holds (for degree $d$ after $M$ prolongations) \cite{vollmer_thesis_2016}:
\begin{equation}\label{eqn:number-unknowns-tomimatsu-sato}
 \begin{split}
  \nvars{d}{M} &= \sum_{l=0}^{\frac{d-e-\tilde{e}}{2}} (2l+1+\tilde{e})(d+1-\tilde{e}-2l)\,\binom{M+3}{2}	\\
  		&= \frac{1}{24}\ (d+2-\Upsigma)\,(M+2)(M+3)\cdot \\
  		&\qquad\qquad\qquad \cdot(d^2+d\Upsigma-2\Upsigma^2+4d+6e(\Upsigma-1)+2\Upsigma+6),
 \end{split}
\end{equation}
where we define $\Upsigma=e+\tilde{e}$.
Hence, if~$d$ is odd, then~$e$ cancels from \eqref{eqn:number-unknowns-tomimatsu-sato}, because $\tilde{e}+e=1$ and thus $\Upsigma=1$ (recall that~$\tilde{e}=\parity(d+e)$).
Similarly, for the number of equations, we find
\begin{equation}\label{eqn:number-equations-tomimatsu-sato}
 \begin{split}
  \meqns{d}{M}	&= \sum_{l=0}^{\frac{d+e-\tilde{e}}{2}} (2l+1+\tilde{e})(d+2-\tilde{e}-2l)\,\binom{M+2}{2}  \\
  		&= \frac{1}{24}\ (d+2+\Updelta)\,(M+1)(M+2)\cdot \\
  		&\qquad\qquad\qquad \cdot(d^2-d\Updelta-2\Updelta^2+6e\Updelta+7d-5\Updelta+12),
 \end{split}
\end{equation}
with $\Updelta=e-\tilde{e}$. Thus, $e$ cancels from \eqref{eqn:number-equations-tomimatsu-sato} if~$d$ is even, because $\tilde{e}=e$ and $\Updelta=0$.

\subsubsection{Kerr metric}\label{sec:kerr}
For the Kerr metric, it is known that a nontrivial Killing tensor of valence~2 exists. We consider the particular case of the extreme Kerr metric (with rotation parameter $a=1$). In Boyer-Lindquist coordinates $(r,\theta,\phi,t)$, this metric reads as follows\cite{boyer_maximal_1967,stephani_exact_2003}:
\begin{equation}\label{eqn:extreme-kerr-metric}
 g = \left(\begin{matrix}
           \frac{r^2+\cos^2(\theta)}{r^2-2r+1} &&& \\
           & r^2+\cos^2(\theta) && \\
           && \frac{P_a(r,\theta)\sin^2(\theta)}{r^2+\cos^2(\theta)}	& \frac{-2\sin^2(\theta)r}{r^2+\cos^2(\theta)} \\
           && \frac{-2\sin^2(\theta)r}{r^2+\cos^2(\theta)} & -\frac{r^2-2r+\cos^2(\theta)}{r^2+\cos^2(\theta)}
     \end{matrix}\right)
\end{equation}
with $P_a(r,\theta)=\cos^2(\theta)r^2+r^4-2\cos^2(\theta)r+\cos^2(\theta)+r^2+2r$.
Algorithm~\ref{alg:sav} yields the following results for the point of reference $(r,\theta)=(2,\nicefrac{\pi}{4})$, \cf\ \cite{vollmer_thesis_2016}.
\begin{center}
\begin{tabular}{|l|c|c|c|c|c|}
  \hline
  Degree						& 1	& 2	& 3	& 4 \\ \hline
  $e=0$ integrals		& 2	& \textbf{5}	& 8	& 14 \\ \hline
  $e=1$ integrals		& 0	& \textbf{0}	& 0	& 0 \\ \hline
  All integrals			& 2	& \textbf{5}	& 8	& 14 \\ \hline
\end{tabular}
\end{center}
Here, $e=0$ and $e=1$ refer again to the branches according to parity in $(p_x,p_y)$.

In second degree we find an upper bound of 5. This is one above the number of trivial Killing tensors of valence~2.
And indeed, the Kerr metric has an additional Killing tensor of valence~2 that commutes with the trivial Killing tensors.
This is the Carter constant \cite{carter_hamilton-jacobi_1968, carter_global_1968}.

\subsubsection{The Darmois solution}\label{sec:darmois}
The Darmois metric is a particular Zipoy-Voorhees metric \cite{stephani_exact_2003,darmois_equations_1927}.
It is therefore a static Tomimatsu-Sato metric.
In prolate spheroidal coordinates it has the following form:
\begin{equation}\label{eqn:zipoy-voorhees}
\resizebox{0.9\textwidth}{!}{
 $g= \left(\frac{x+1}{x-1}\right)^2
          \left((x^2-y^2)\left(\frac{x^2-1}{x^2-y^2}\right)^4
                  \left(\frac{dx^2}{x^2-1}+\frac{dy^2}{1-y^2}\right)
                +(x^2-1)(1-y^2)d\phi^2\right)
        -\left(\frac{x-1}{x+1}\right)^2\ dt^2$.
}
\end{equation}

Existence of a nontrivial Killing tensor for the Darmois solution has been suggested by \cite{brink_II_2008, brink_formal_2009}, but later studies challenged this claim \cite{kruglikov_nonexistence_2011,gerakopoulos_non-integrability_2012,maciejewski_nonexistence_2013}.
In~\cite{kruglikov_nonexistence_2011}, Kruglikov \& Matveev study the number of Killing tensors for this metric up to valence~6, and the method we discuss here is based on this work. It is therefore interesting to compare the computer performance of both methods (see below).

For static and axially symmetric metrics, the algorithm can be improved further.
It is possible to restrict to integrals of even parity in $p_\phi$. Their parity \wrt\ $(p_x,p_y)$ can be taken equal to the parity of $d$. In order to make a statement for valence~$d$, we need to apply the algorithm for $d$, $d-1$ and $d-2$.

The reason for this simplification is the following observation:
The Hamiltonian is not only of even parity \wrt\ $(p_x,p_y)$, but also \wrt\ $p_\phi$ or~$p_t$. Thus, components of the integral of even and odd parity \wrt\ $p_\phi$ can be considered separately.

Now, consider only integrals of even parity in $p_\phi$ and such that their parity \wrt\ $(p_x,p_y)$ equals the parity of $d$ as an integer number.
Let~\Sd\ be the system of equations obtained from the Poisson equation by considering coefficients with respect to momenta.
For Weyl metrics, the system of equations \Sd\ splits into four separate subsystems, see Figure~\ref{fig:subsystems}. Each of the subsystems can be solved independently.
Unknowns of one of the subsystems do not appear in the other subsystems.
The crucial observation is that each of this subsystems corresponds to one of the described type, but possibly with another value for $d$. For details, see~\cite{vollmer_thesis_2016}.

\begin{figure}
\begin{center}
\textbf{Subsystems according to parity of the integral}
\smallskip

\begin{tabular}{|c|c|c|}
\hline
\textbf{parity \wrt\ $(p_x,p_y)$}	& \textbf{even parity in $p_\phi$}	& \textbf{odd parity in $p_t$}	\\ \hline
{equals $\parity(d)$}							& \Sd																& corresponds to $\Seqn_{d-2}$	\\ \hline
{opposite to $\parity(d)$}				& corresponds to $\Seqn_{d-1}$			& corresponds to $\Seqn_{d-1}$	\\ \hline
\end{tabular}	
\end{center}
\caption{Refined parity split for the Darmois metric (and the C-metric). For details, see~\cite{vollmer_thesis_2016}.}
\label{fig:subsystems}
\end{figure}

\begin{theorem}\label{thm:ZV-11}
The Darmois solution, \ie\ the Zipoy-Voorhees space-time with parameter $\delta=2$, has no additional independent Killing tensor of valence~$d\leq 11$ that is in involution with the trivial Killing tensors $d\phi$, $dt$, and the metric.
\end{theorem}

\begin{center}
\textbf{Results Darmois metric}\\
\textit{(all values refer to \Sd\ only)}
\vskip 0.2cm
\begin{tabular}{|c|c|cccccc|}
\hline
d   & \upperbound[d]{d} & $\meqns{d}{d}$	& $\nvars{d}{d}$	& rows of $M$ & columns of $M$ 	& $\rg(M)$	& time  \\ \hline
9		&	0									& 5005						& 4620						& 1058				& 726							&	726				& 1.8h	\\ \hline
10	& 21								& 7392						& 7098						& 1358				& 1064						& 1043			& 13h		\\ \hline
11	& 0									& 10780						& 10192						&	2162				& 1510						& 1510			& 7days	\\ \hline
\end{tabular}
\end{center}
The upper bound to the number of Killing tensors is therefore 42, and this equals the number of minimally expected Killing tensors (note that $\Seqn_{d-1}$ contributes twice).
The point of reference in this computation is $(\nicefrac12,2)$.

\paragraph{Computational Efficiency}\label{par:efficiency}
Let us quickly contrast the performance of the method presented in Section~\ref{sec:methods} with the original method from \cite{kruglikov_nonexistence_2011}.

While on our computer the algorithm from~\cite{kruglikov_nonexistence_2011} takes about 3 minutes for all the necessary computations for valence $d=4$, the modified algorithm takes only around 3 seconds.
For valence $d=6$, we needed approximately 47 minutes with the algorithm from~\cite{kruglikov_nonexistence_2011}. We needed a little more that a minute (67 seconds) with the improved method.
These numbers include both branches for~\cite{kruglikov_nonexistence_2011} and all three branches for the method of Algorithm~\ref{alg:sav} in combination with the decomposition according to Figure~\ref{fig:subsystems}.


\subsection{C-metrics}

The C-metric is used in relativity as a model for systems of two black holes accelerating in opposite directions under the action of certain forces.
Its physical properties are described, for instance, in~\cite{griffiths_interpreting_2006}.
We express the C-metric in a form given by Hong and Teo \cite{hong_form_2003,griffiths_interpreting_2006},
\begin{equation}\label{eqn:c-metric}
 g = \frac{1}{\alpha^2(x+y)^2}\,\left(\frac{dx^2}{X}+\frac{dy^2}{Y}+X\,d\phi^2-Y\,d\tau^2\right),
\end{equation}
where $X$ and $Y$ are
\[
 X=(1-x^2)(1+2m\alpha\,x), \qquad Y=(y^2-1)(1-2m\alpha\,y).
\]
This metric can describe space-times with different properties \cite{griffiths_interpreting_2006}.
The most interesting case is the case of two accelerated black holes and requires $x+y>0$ and $-1<x<1$ \cite{griffiths_interpreting_2006}.
The parameter $\alpha$ describes an acceleration and $m$ denotes a mass parameter; they are restricted by the requirement $0<2m\alpha<1$ \cite{griffiths_interpreting_2006}.
Space-times of the kind described have three regions, separated by what is called the black hole horizon (at $y=1$) and the acceleration horizon (at $y=\frac{1}{2m\alpha}$).
In the region between the horizons, the space-time is static and $\partial_t$ is a timelike Killing vector \cite{griffiths_interpreting_2006}.
We only consider this region for brevity of discussion.

\begin{theorem}\label{thm:C-metric}
The C-metric with $\alpha=\nicefrac12$ and $m=\nicefrac12$ has no additional, independent Killing tensor of valence~$d\leq 9$ that is in involution with $d\phi$, $dt$, and the metric.
\end{theorem}

\begin{center}
\textbf{Results C-metric with $\alpha=\nicefrac12$ and $m=\nicefrac12$}\\
\textit{(all values refer to \Sd\ only)}
\vskip 0.2cm
\begin{tabular}{|c|c|cccccc|}
\hline
d   & \upperbound[d]{d} & $\meqns{d}{d}$	& $\nvars{d}{d}$	& rows of $M$ & columns of $M$ 	& $\rg(M)$	& time  \\ \hline
7		&	0									& 1980						& 1800						& 488					& 308							&	308				& 17m		\\ \hline
8		& 15								& 3150						& 3025						& 608					& 468							& 468				& 21h		\\ \hline
9		& 0									& 5005						& 4620						&	1113				& 728							& 728				& 57h		\\ \hline
\end{tabular}
\end{center}
The upper bound to the number of Killing tensors is therefore 30, and this equals the number of minimally expected Killing tensors.
The point of reference for the computation was $(0,\nicefrac32)$.

\section{Conclusions}
We have seen that taking into account the full structural properties of the system of PDEs connected to~\eqref{eqn:killing}, the method suggested by \cite{kruglikov_nonexistence_2011} can be modified to achieve higher computational efficiency. In particular, we introduced an elimination scheme and the freedom to add trivial integrals to speed up the computations and to reach higher degrees of the integrals. Moreover, we saw that for static space-times the Killing tensors partially come from lower degrees, and that this fact can significantly improve computational performance.

The additional techniques were implemented for several examples, namely the Darmois solution, a C-metric and a Tomimatsu-Sato metric. The extreme Kerr solution was discussed as an example that admits a nontrivial Killing tensor.
We saw that the new techniques also work fine with axially symmetric metrics that are stationary (instead of only static).

Obviously, the techniques applied in our method rely on the structural properties of the problem rather than its physical properties.
The method is applicable not only for stationary and axially symmetric metrics, and it can produce results in other contexts as well, see~\cite{kruglikov_certain_2015} for an application of a similar method in sub-Riemannian geometry.

\section*{Acknowledgements}
The author wishes to thank Vladimir Matveev for introduction to the pro\-lon\-gation-projection approach to the existence of Killing tensors, as well as for generous support and fruitful discussions. He is grateful for financial support by the Research Training Group Quantum and Gravitational Fields at the University of Jena, financed by Deutsche Forschungsgemeinschaft.
The author wishes to thank Galliano Valent for drawing his attention to the C-metric.

\sloppy
\printbibliography

\end{document}